\documentclass[12pt]{article}
\usepackage[latin9]{inputenc}
\usepackage[letterpaper]{geometry}
\usepackage{float}
\usepackage{amsmath}
\usepackage{amsthm}
\usepackage{graphicx}

\makeatletter
\numberwithin{equation}{section}
\numberwithin{figure}{section}
\theoremstyle{plain}
\newtheorem*{thm*}{\protect\theoremname}

\@ifundefined{date}{}{\date{}}

\theoremstyle{theorem}

\theoremstyle{definition}

\makeatother

\providecommand{\theoremname}{Theorem}

\begin{document}
\begin{center}
{\large{}The altitudes of a triangle }{\large\par}
\par\end{center}

\begin{center}
Mark Mandelkern
\par\end{center}

\vspace{0.3cm}

\begin{quotation}
\noindent \emph{Abstract. }A long-standing, unanswered question regarding
Euclid's \emph{Elements} concerns the absence of a theorem for the
concurrence of the altitudes of a triangle, and the possible reasons
for this omission. In the centuries following Euclid, a remarkable
number of proofs have been put forward; this suggests a search for
the most elementary and direct proof. This paper provides a simple,
direct, elementary proof of the theorem; it is based solely on the
\emph{Elements}. 
\end{quotation}
\vspace{0.8cm}

The concurrence of the angle bisectors of a triangle, at the \emph{incenter,}
follows directly from Proposition 4 in Book IV of Euclid's \emph{Elements
}{[}2{]}, and the concurrence of the perpendicular bisectors of the
sides, at the \emph{circumcenter,} follows directly from Proposition
5. In the extant Euclid, however, there is no mention anywhere of
the concurrence of the altitudes, at the \emph{orthocenter.} This
omission has long been a source of mystery and speculation. In the
centuries following Euclid, a remarkable number of proofs have been
put forward; this suggests a search for the most elementary and direct
proof. 

A history of the altitude theorem, a comprehensive discussion of known
proofs, and an extensive list of references will be found in {[}1{]}.\emph{
}Several proofs listed there involve the drawing of additional lines,
some depend on a concurrence theorem for the incenter or circumcenter,
and others depend on Ceva's Theorem, vectors, norms, complex numbers,
trigonometry, analytical geometry, Cartesian coördinates, and other
methods introduced only in recent centuries. Only two of the listed
proofs could be considered elementary. Proof 3, based on the theory
of cyclic quadrilaterals, is not strictly direct, nor completely elementary.
The most elementary of the listed proofs is Proof 2, based on the
theory of similar triangles; it is attributed to Newton. 

The proof here, while also involving similar triangles, differs from
Newton's proof in a fundamental manner. The proof attributed to Newton
notes the intersection points of two different pairs of altitudes,
and shows that these  points coincide. The proof here uses the strategy
(also used in {[}1, Proof 3{]}) of drawing only two of the altitudes,
and then proving that the line drawn from the third vertex, through
the point where these two meet, is the third altitude; this technique
is more direct and elementary.

The proof  below is simple, direct, and elementary; it introduces
no additional lines, depends on none of the advanced methods of modern
times, respects the axiomatic character of the \emph{Elements}, and
could have been included by Euclid. 

\noindent {\footnotesize{}\rule[0.5ex]{0.45\columnwidth}{0.5pt}}{\footnotesize\par}

\noindent {\footnotesize{}2020 Mathematics Subject Classification.
Primary 51M04. }{\footnotesize\par}

\noindent {\footnotesize{}Key words and phrases. Triangle, altitude,
concurrence. }{\footnotesize\par}

\begin{figure}[H]
\includegraphics[viewport=100bp 130bp 700bp 600bp,scale=0.4]{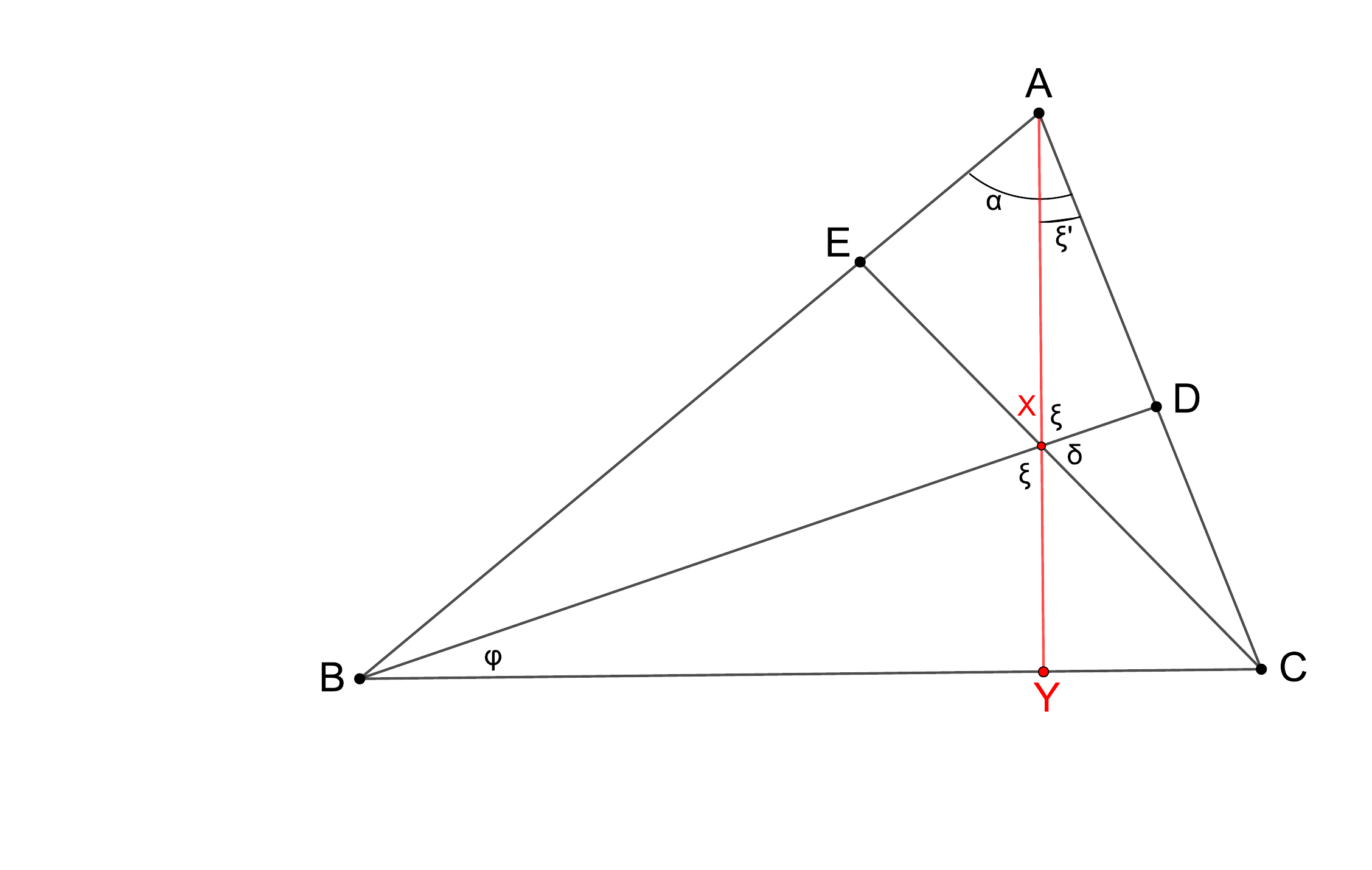}
\end{figure}

\begin{thm*}
The altitudes of a triangle are concurrent. 
\end{thm*}
\begin{proof}
Let $ABC$ be an acute-angled triangle. (The proof for an obtuse triangle
will be essentially the same, \emph{mutatis mutandis. }Alternatively,
the Lemma in {[}1{]} could be used to reduce the problem to acute
triangles.) 

Draw the altitudes $BD$ and $CE$, intersecting at the point $X$,
and draw the line $AY$ through $X$ to the point $Y$ on side $BC$.
Write $\varphi=\angle CBD$, and $\xi=\angle BXY$; clearly, $\angle AXD=\xi$,
and $\angle XAD=\xi'$, the complement of $\xi$. The angles $\alpha=\angle CAE$
and $\delta=\angle CXD$ are each complementary to angle $\angle ACE$;
thus $\alpha=\delta$. 

The right triangles $ADB$ and $XDC$ have the corresponding equal
acute angles $\alpha,\delta$; thus they are equal-angled. The proportion\footnote{In keeping with the character of the \emph{Elements}, this proportion
is written in a manner which avoids reference to algebra in the modern
sense.} $AD:BD::XD:CD$ follows from {[}VI,4{]}.\footnote{This notation refers to Book VI, Proposition 4, in Euclid's \emph{Elements
}{[}2{]}.} Noticing that the alternate proportion $AD:XD::BD:CD$, derived from
{[}V,16{]}, expresses a proportionality for the sides of the right
triangles $ADX$ and $BDC$ , and applying {[}VI,6{]}, we find that
these triangles are equal-angled; thus $\varphi=\xi'$. 

Now two of the angles in triangle $BYX$ are complementary; thus the
third angle $\angle BYX$ is a right angle. This shows that the line
$AY$ is the third altitude. \\
\end{proof}
\noindent References

\noindent {\small{}{[}1{]} Hajja, M., Martini, H., Concurrency of
the altitudes of a triangle, }\emph{\small{}Math. Semesterber. }{\small{}60(2013),
249--260. }{\small\par}

\noindent {\small{}{[}2{]} Heath, T. L., }\emph{\small{}The Thirteen
Books of Euclid's Elements, }{\small{}2nd ed., 3 volumes, Cambridge
University Press, 1926. Reprint, Dover, 1956. Reprint, Cambridge University
Press, 2015. }\\

\noindent {\small{}New Mexico State University, Las Cruces, New Mexico,
USA }{\small\par}

\noindent {\small{}http://www.zianet.com/mandelkern}{\small\par}

\noindent {\small{}mandelkern@zianet.com, mmandelk@nmsu.edu }\\

\end{document}